\documentclass[12pt,thmsa]{article}
\usepackage{amsmath, latexsym, amsfonts, amssymb, amsthm, amscd}

\newtheorem{theorem}{Theorem}
\newtheorem{corollary}{Corollary}
\newtheorem{proposition}{Proposition}
\newtheorem{lemma}{Lemma}

\newcommand{\p}{\Bbb{P}}
\newcommand{\px}{\Bbb{P}_x}
\newcommand{\e}{\Bbb{E}}

\newcommand{\ind}{\mbox{\rm 1\hspace{-0.04in}I}}

\newcommand{\R}{\mbox{\rm I\hspace{-0.02in}R}}

\newcommand{\ud}{\textrm{d}}

\def\QED{\hfill\vrule height 1.5ex width 1.4ex depth -.1ex \vskip20pt}

\begin{document}
\hspace*{-0.5in}{\footnotesize This version August 17, 2007.}
\vspace*{0.9in}
\begin{center}
{\Large Some explicit identities associated with positive\\
self-similar Markov processes.}\\
\vspace*{0.4in} {\large L. Chaumont\footnote{ LAREMA, D\'epartement
de Math\'ematiques, Universit\'e d'Angers. 2, Bd Lavoisier - 49045,
{\sc Angers Cedex 01. France.} E-mail:
loic.chaumont@univ-angers.fr.}, A.E. Kyprianou\footnote{$^{,3}$
Department of Mathematical Science, University of Bath. {\sc Bath,
BA2 7AY. United Kingdom.}
$^{2}$E-mail: ak257@bath.ac.uk.} and J.C.
Pardo$^{*,}$}\footnote{Laboratoire de Probabilit\'es et Mod\`eles Al\'eatoires, Universit\'e Pierre et Marie Curie, 4, Place Jussieu - 75252 PARIS CEDEX 05. E-mail: jcpm20@bath.ac.uk.\\
Reasearch supported by EPSRC grant EP/D045460/1.\\
$^*$Corresponding author.}
\end{center}
\vspace{0.2in}

\begin{abstract} We consider some special classes of L\'evy processes with no gaussian component
whose L\'evy measure is of the type $\pi(dx)=e^{\gamma
x}\nu(e^x-1)\,dx$, where $\nu$ is the density of the stable L\'evy
measure and $\gamma$ is a positive parameter which depends on its
characteristics. These processes were introduced in \cite{CC} as the
underlying L\'evy processes in the Lamperti representation of
conditioned stable L\'evy processes. In this paper, we compute
explicitly the law of these L\'evy processes at their first exit
time from a finite or semi-finite interval, the law of their
exponential functional and the
first hitting time probability of a pair of points.\\

\noindent {\sc Key words and phrases}: Positive self-similar Markov
processes, Lamperti representation, conditioned stable L\'evy
processes, first exit time, first hitting time, exponential functional.\\

\noindent MSC 2000 subject classifications: 60 G 18, 60 G 51, 60 B
52.
\end{abstract}

\vspace{0.5cm}

\section{Introduction}
In recent years there has been a general recognition that L\'evy
processes play an ever more important role in various domains of
applied probability theory such as financial mathematics, insurance
risk, queueing theory, statistical physics  or mathematical biology.
In many instances there is a need for explicit examples of L\'evy
processes where tractable mathematical expressions in terms of the
characteristics of the underlying L\'evy process may be used for the
purpose of numerical simulation. Depending on the problem at hand, particular
functionals are involved such as the first entrance
times and overshoot distributions.

In this paper, we exhibit some special classes of L\'evy processes
for which we can compute explicitly the law of the position at the
first exit time of an interval, the two points hitting probability
and the exponential functional. Moreover, two new, concrete examples of
 scale functions for  spectrally one sided processes will fall out of
our analysis.

Known examples of overshoot distributions concern essentially (some particular
classes) of strictly stable processes and processes whose jumps are
of a compound Poisson nature with exponential jumps (or slightly
more generally whose jump distribution has a rational Fourier
transform). For example,  let us state the solution of the two sided
exit problem for completely asymmetric stable processes. In that
case we take $(X,\mathbf{P}_x)$,  $x\in \mathbb{R}$, to be a
spectrally positive L\'evy stable process with index
$\alpha\in(1,2)$ starting from $x$.  Let $\sigma^+_a = \inf\{t>0 :
X_t
>a \}$ and $\sigma^-_0 = \inf\{t>0 : X_t <0\}$. It is known (cf.
Rogozin \cite{ro}) that for $y>0$,
\[
\mathbf{P}_x \Big(X_{\sigma^+_a}-a\in \ud y; \sigma^+_a <
\sigma^-_0\Big) = \frac{\sin\pi(\alpha-1)} {\pi}\left(\frac{a-x}{x
y}\right)^{(\alpha-1)} \frac{\ud y}{(y+a)(y+a-x)}.
\]
For the case of processes whose jumps are of a compound Poisson
nature with exponential jumps, the overshoot distribution is
again exponentially distributed; see Kou and Wang \cite{KW}.
See also Lewis and Mordecki \cite{LM} and Pistorius \cite{P} for the
more general case of a jump distribution with a rational Fourier
transform and for which the overshoot distribution belongs to the same
class as the respective jump distribution of the underlying L\'evy
process.

The exponential functional of a L\'evy process, $\xi$, i.e.
\[
\int_0^\infty\exp\Big\{-\xi_s\Big\}\,\ud s ,
\]
also appears in various aspects of probability theory, such as: 
self-similar Markov processes, random processes in random
environment, fragmentation processes, mathematical finance, Brownian
motion on hyperbolic spaces, to name but a few. 
 In general, the distribution of exponential functionals can be rather
complicated. Nonetheless, it is known for the case that $\xi$ is either: a standard Poisson processes,  Brownian motion with drift
 and a particular class of spectrally
negative L\'evy processes of bounded variation whose Laplace
exponent is of the form
\[
\psi(q)=\frac{q(q+1-a)}{b+q}, \qquad q\geq 0,
\]
where $0<a<1<a+b$. See Bertoin and Yor \cite{BY} for an overview 
on this topic.

The class of L\'evy processes that we consider in this paper do not fulfill a scaling
property and may have two-sided jumps. Moreover, they have no
Gaussian component and their L\'evy measure is of the type
$\pi(dx)=e^{\gamma x}\nu(e^x-1)\,dx$, where $\nu$ is the density of
the stable L\'evy measure with index $\alpha\in(0,2)$ and $\gamma$
is a positive parameter which depends on its characteristics.  It is
not difficult to see that the latter L\'evy measure has a density
which is asymptotically equivalent to that of an $\alpha$-stable
process for small $|x|$ and has exponential decay for large $|x|$.
This implies that such processes have paths which are of bounded or
unbounded variation accordingly as $\alpha\in(0,1)$ and $\alpha\in
[1,2)$ respectively. Further, they also have exponential moments.
Special families of tempered stable processes, also known as CGMY
processes, are classes of L\'evy processes with similar properties
to the aforementioned which have enjoyed much exposure in the
mathematical finance literature as instruments for modelling risky
assets. See for example Carr et al. \cite{CGMY}, Boyarchenko and
Levendorskii \cite{boy}, Cont \cite{cont} or Schoutens \cite{Sch}.
Although the L\'evy processes presented in this paper are not
tempered stable processes, it is intriguing to note that they
possess properties which have proved to be popular for financial
models but now with the additional luxury that they come with a
number of explicit fluctuation identities.

We conclude the introduction with a brief outline of the remainder
of the paper. The next section introduces the classes of processes
which are concerned in this study. In section 3, we give the law of
the position at the first exit time from a (semi-finite) interval.
In section 4 we compute explicitly the two point hitting probability
and in section 5, we study the law of the exponential functional of
L\'evy-Lamperti processes.

\section{Preliminaries on L\'evy-Lamperti processes}\label{prelim}

Denote by $\mathcal{D}$ the Skorokhod space of $\mathbb{R}$-valued
c\`adl\`ag paths and by $X$ the canonical process of the coordinates
on $\mathcal{D}$. Positive ($\mathbb{R}_+$-valued), self-similar
Markov processes $(X,\p_x)$, $x>0$, are strong Markov processes with
paths in $\mathcal{D}$, which fulfill a scaling property, i.e.~there
exists a constant $\alpha > 0$ such that for any $b>0$:
\begin{equation}\label{scale}
\mbox{\it The law of $\;(bX_{b^{-\alpha}t},\,t\ge0)$ under $\p_x$ is
$\p_{bx}$.}
\end{equation}
We shall refer to these processes  as pssMp. According to Lamperti
\cite{La}, any pssMp up to its first hitting time of 0 may be
expressed as the exponential of a L\'evy process, time changed by
the inverse of its exponential functional. More formally, let
$(X,\p_x)$ be a pssMp with index $\alpha>0$, starting from $x>0$,
set
\[S = \inf \{t>0 : X_t =0\}\]
and write the canonical process $X$ in the following form:
\begin{equation}\label{lamp}
X_t=x\exp\left\{\xi_{\tau(tx^{-\alpha})}\right\} \qquad 0\le t<S\,,
\end{equation}
where for $t<S$,
\[\tau(t) = \inf \left\{s\geq 0 : \int_0^s
\exp\left\{\alpha\xi_u\right\} \ud u \geq t\right\}.\] Then under
$\p_x$, $\xi=(\xi_t,\;t\geq 0)$ is a L\'evy process started from $0$
whose law does not depend on $x>0$ and such that:
\begin{itemize}
\item[$(i)$] if  $\p_x(S=+\infty)=1$,  then  $\xi$ has an infinite lifetime
and $\displaystyle\limsup_{t\rightarrow+\infty}\xi_t=+\infty$,
$\p_x$-a.s.,

\item[$(ii)$] if $\p_x(S<+\infty,\,X(S-)=0)=1$, then $\xi$ has an infinite
lifetime and $\displaystyle\lim_{t\to\infty} \xi_t = -\infty$,
$\p_x$-a.s.,

\item[$(iii)$] if $\p_x(S<+\infty,\,X(S-)>0)=1$,  then $\xi$ is killed at an
independent exponentially distributed random time with parameter
$\lambda>0$.

\end{itemize}
As it is mentioned in \cite{La}, the probabilities
$\p_x(S=+\infty)$, $\p_x(S<+\infty,\,X(S-)=0)$ and
$\p_x(S<+\infty,\,X(S-)>0)$ are 0 or 1 independently of $x$, so that
the three classes presented above are exhaustive. Moreover, for any
$t < \int_0^{\infty}\exp\{\alpha\xi_s\}\,\ud s$,
\begin{equation}\label{1664}\tau(t)=\int_0^{x^\alpha t}\frac{\ud s}
{(X_s)^\alpha}\,,\;\;\; \p_x-\mbox{a.s.}\end{equation} Therefore
(\ref{lamp}) is invertible and yields a one to one relation between
the class of pssMp's killed at time $S$ and the one of L\'evy
processes.\\

Now let us consider three particular classes of pssMp. (We refer to
\cite{CC} for more details in what follows.) The first one is
identified as a  stable L\'evy processes killed when it first exits from
the positive half-line. In particular, if $\mathbf{P}_x$ is the law
of a stable L\'evy process with index $\alpha$ (or $\alpha$-stable
process for short) initiated from $x>0$ with $\alpha\in(0,2]$, then
with $T=\inf\{t:X_t\le0\}$, under $\mathbf{P}_x$, the process
\[X_t\ind_{\{t<T\}}\]
is a pssMp which satisfies condition $(ii)$ if it has no negative
jumps or $(iii)$ if it has negative jumps. We call $\xi^*$ the
L\'evy process (with finite or infinite lifetime) resulting from the
Lamperti representation of the killed stable process. The
characteristic exponent of $\xi^*$ has been computed in \cite{CC}
and is given by
\begin{equation}\label{phistar}
\Phi^*(\lambda)=ia^*\lambda+\int_{\mathbb{R}}[e^{i\lambda
x}-1-i\lambda(e^x-1)\ind_{\{|e^x-1|<1\}}]\pi^*(x)\,\ud
x-c_{-}\alpha^{-1}\,, \;\;\;\lambda\in\mathbb{R}\,,\end{equation}
where $a^*$ is a constant,
\[\pi^*(x)=\frac{c_+e^x}{(e^x-1)^{\alpha+1}}\ind_{\{x>0\}}+
\frac{c_-e^x}{(1-e^x)^{\alpha+1}}\ind_{\{x<0\}}\,,\]
 and $c_-$, $c_+$ are nonnegative
constants such that $c_-c_+>0$. Note that the L\'evy measure of
$\xi^*$ satisfies $\pi^*(x)=e^x\nu(e^x-1)$, where $\nu$ is the
density of the stable L\'evy measure with index $\alpha$ and
symmetry parameters $c_-$ and $c_+$.

The second class is that of stable processes conditioned
to stay positive. ( See for instance
in \cite{Ch} for an overview of such processes.) A process in this class is the result of a Doob $h$-transform with $h(x) =
x^{\alpha\rho}$ and $\rho=\mathbf{P}_0(X_1<0)$. More precisely, $h$
is invariant for the killed process mentioned above
$(X_t\ind_{\{t<T\}},\mathbf{P}_x)$ and the law $\p^\uparrow_x$
defined on each $\sigma$-field $\mathcal{F}_t$ generated by the
canonical process up to time $t$ by
\begin{equation}\label{uparrow}
\left.\frac{\ud\mathbb{P}^\uparrow_x}{\ud\mathbf{P}_x}\right|_{\mathcal{F}_t}
= \frac{X^{\alpha\rho}_t}{x^{\alpha\rho}}\mathbf{1}_{\{t<T\}}
 \end{equation}
is this of a pssMp which derives toward $+\infty$ (in particular it
satisfies condition $(i)$). Then the underlying L\'evy process,
which will be denoted by $\xi^\uparrow$, is such that
\[
\lim_{t\rightarrow+\infty}\xi^\uparrow_t=+\infty,\qquad
\textrm{a.s.,}
\]
and from \cite{CC} its characteristic exponent is
\begin{equation}\label{phiup}
\Phi^\uparrow(\lambda)=ia^\uparrow\lambda+\int_{\mathbb{R}}[e^{i\lambda
x}-1-i\lambda(e^x-1)\ind_{\{|e^x-1|<1\}}]\pi^\uparrow(x)\,\ud x\,,
\;\;\;\lambda\in\mathbb{R}\,,\end{equation} where $a^\uparrow$ is a
real constant and
\[\pi^\uparrow(x)=\frac{c_+e^{(\alpha\rho+1)x}}{(e^x-1)^{\alpha+1}}\ind_{\{x>0\}}+
\frac{c_-e^{(\alpha\rho+1)x}}{(1-e^x)^{\alpha+1}}\ind_{\{x<0\}}\,.\]

The third class of pssMp that we will consider is that of stable
processes conditioned to hit 0 continuously. Processes in this class are again it is defined as
a Doob $h$-transform with respect to the function $h'(x)=\alpha\rho
x^{\alpha\rho-1}$ which is also invariant for the killed process
$(X_t\ind_{\{t<T\}},\mathbf{P}_x)$. Then  the law $\p^\downarrow_x$
which is defined on each $\sigma$-field $\mathcal{F}_t$ by
\begin{equation}\label{downarrow}
\left.\frac{\ud\mathbb{P}^\downarrow_x}{\ud\mathbf{P}_x}\right|_{\mathcal{F}_t}
=
 \frac{X^{\alpha\rho-1}_t}{x^{\alpha\rho-1}}\mathbf{1}_{\{t<T\}}
\end{equation}
 is this of a pssMp who hits 0 in a continuous way, i.e.
$(X,\p_x^\downarrow)$ satisfies condition $(ii)$. Let
$\xi^\downarrow$ by the underlying L\'evy process in the Lamperti
representation of this process, then
\[
\lim_{t\rightarrow+\infty}\xi^\downarrow_t=-\infty\qquad \textrm{a.s.}, 
\]
and the
characteristic exponent of $\xi^\downarrow$ is given by
\begin{equation}\label{phiup}
\Phi^\downarrow(\lambda)=ia^\downarrow\lambda+\int_{\mathbb{R}}[e^{i\lambda
x}-1-i\lambda(e^x-1)\ind_{\{|e^x-1|<1\}}]\pi^\downarrow(x)\,\ud x\,,
\;\;\;\lambda\in\mathbb{R}\,,\end{equation} where $a^\downarrow$ is
a constant and
\[\pi^\downarrow(x)=\frac{c_+e^{\alpha\rho x}}{(e^x-1)^{\alpha+1}}\ind_{\{x>0\}}+
\frac{c_-e^{\alpha\rho x}}{(1-e^x)^{\alpha+1}}\ind_{\{x<0\}}\,.\]
Note that the constants $a^*$, $a^\uparrow$ and $a^\downarrow$ are
computed explicitly in \cite{CC} in terms of $\alpha$, $\rho$, $c_-$
and $c_+$. Actually the process $\xi^\downarrow$ corresponds to
$\xi^\uparrow$ conditioned to drift toward $-\infty$ (or
equivalently $\xi^\uparrow$ is $\xi^\downarrow$ conditioned to drift
to $+\infty$). We will sometime use this relationship which is
stated in a more formal way the next proposition. In the sequel, $P$
will be a reference probability measure on ${\cal D}$ under which
$\xi^*$, $\xi^\uparrow$ and $\xi^\downarrow $ are L\'evy processes
whose respective laws are defined above.
\begin{proposition}\label{4590}
For every $t\ge0$, and every bounded measurable function $f$,
\[E[f(\xi^\uparrow_t)]=
E[\exp({\xi_t^\downarrow})f(\xi^\downarrow_t)]\,.\] In particular,
processes $-\xi^\uparrow$ and  $\xi^\downarrow$ satisfy Cramer's
condition: $E(\exp{-\xi_1^\uparrow})=1$ and
$E(\exp{\xi_1^\downarrow})=1$.
\end{proposition}
\noindent{\it Proof}. Let $f$ be as in the statement. From
(\ref{uparrow}) and (\ref{downarrow}), we deduce that for every
$\mathbb{P}^\downarrow_x$-a.s. finite $({\cal F}_u)$-stopping time
$U$,
\begin{equation}\label{3521}
x\mathbb{E}^\uparrow_x[f(X_U)]=
\mathbb{E}^\downarrow_x[X_Uf(X_U)]\,.\end{equation} Let $t\ge0$. By
applying (\ref{3521}) to the $({\cal F}_u)$-stopping time
\[x^{\alpha}\inf\Big\{u:\tau(u)>t\Big\},\]
which is $\mathbb{P}^\downarrow_x$-a.s. finite, and using
(\ref{lamp}) (note that $\tau(u)$ is continuous and increasing), we
obtain
\[E[f(\xi^\uparrow_t)]=
E[\exp(\xi_t^\downarrow)f(\xi^\downarrow_t)]\,,\]  which is the
desired result. \QED \noindent We refer to Rivero \cite{ri}, IV.6.1
for a similar discussion on conditioned stable processes considered
as pssMp. In the sequel we call $\xi^*$, $\xi^\uparrow$ and
$\xi^\downarrow$ the {\it L\'evy-Lamperti processes}. We now compute
the law of some of their functionals.

\section{Entrance laws for L\'evy-Lamperti processes: intervals}

In this section, by studying the two sided exit problems for
$\xi^\uparrow$, $\xi^*$ and $\xi^\downarrow$, we shall obtain a
variety of new identities including the identification of two new
scale functions in the case of one-sided jumps.

To this end, we shall start with a generic result pertaining to any
positive self-similar Markov process  $(X,\mathbb{P}_x)$, for $x>0$. Recall that $P$ is
the reference probability measure on $D$. Let $\xi$ be a L\'evy process starting from $0$, under $P$, with the
same law as the underlying L\'evy process associated to   $(X,\mathbb{P}_x)$. For any $y\in\mathbb{R}$ let
\[
T^+_y=\inf\{t:\xi_t\ge y\}\;\;\mbox{and}\;\;T_y^-=\inf\{t:\xi_t\le
y\}\,,
\]
and for any $y>0$ let
\[
 \sigma^+_y=\inf\{t:X_t\ge
y\}\;\;\mbox{and}\;\;\sigma_y^-=\inf\{t:X_t\le y\}.
\]

\begin{lemma}\label{generic}Fix $-\infty< v<0<u<\infty$.
Suppose that $A$ is any interval in $[u,\infty)$ and $B$ is any
interval in $(-\infty, v]$. Then,
\[
 P\Big(\xi_{T^+_u}\in A; T^+_u< T^-_v\Big) =
 \mathbb{P}_1\Big(X_{\sigma^+_{e^u}} \in e^A; \sigma^+_{e^u}< \sigma^-_{e^v}\Big)
\]
and
\[
P\Big(\xi_{T^-_v}\in B; T^+_u> T^-_v\Big) =
\mathbb{P}_1\Big(X_{\sigma^-_{e^v}} \in e^B; \sigma^+_{e^u}>
\sigma^-_{e^v}\Big).
\]
\end{lemma}
The proof is a straightforward consequence of the Lamperti
representation (\ref{lamp}) and is left as an exercise. Although
somewhat obvious, this lemma indicates that for the three processes
$\xi^\uparrow$, $\xi^*$ and $\xi^\downarrow$, we need to understand
how, respectively, an  $\alpha$-stable  process conditioned to stay positive, an
$\alpha$-stable  process killed when it exits  the positive
half-line and an $\alpha$-stable process conditioned to hit the
origin continuously, exit a positive interval around $x>0$.
 Fortunately this is possible thanks to a result of
Rogozin \cite{ro} who established the following result for
$\alpha$-stable processes.
\begin{theorem}[Rogozin \cite{ro}]\label{rogozin}
Suppose that $(X,\mathbf{P}_x)$ is an $\alpha$-stable process, initiated from $x$, which
has two sided jumps. Denoting $\rho
=\mathbf{P}_0(X_1 <0)$ we have for $a>0$ and $x\in (0,a)$,
\begin{eqnarray*}
 &&\hspace{-1cm}\mathbf{P}_x\Big(X_{\sigma^+_a} - a \in \ud y; \sigma^+_a < \sigma^-_0\Big) \\
&&\hspace{1cm}= \frac{\sin\pi\alpha(1-\rho)}{\pi} (a-
x)^{\alpha(1-\rho)}x^{\alpha\rho}y^{-\alpha(1-\rho)}(y+a)^{-\alpha\rho}(y+a-x)^{-1}\ud
y
\end{eqnarray*}
\end{theorem}
\noindent Note that an expression for $ \mathbf{P}_x(-X_{\sigma^-_0}
\in \ud y; \sigma^+_a > \sigma^-_0)$ can be derived from the above
expression by replacing $x$ by $a-x$ and $\rho$ by $1-\rho$.

\bigskip

In the sequel, with an abuse of notation, we will denote by $T^+_y$
and $T^-_y$ for the first passage times above and below $y\in \R$,
respectively, of the processes $\xi^\uparrow, \xi^*$ or
$\xi^\downarrow$ depending on the case that we are studying.

We now proceed to split the remainder of this section into three
subsections dealing with the two sided exit problem and its
ramifications   for the three processes $\xi^\uparrow, \xi^*$ and
$\xi^\downarrow$ respectively.


\subsection{Calculations for $\xi^\uparrow$}\label{xiup}

 The two sided exit problem for $\xi^\uparrow$
 can be obtained from Lemma \ref{generic} and Theorem \ref{rogozin} as follows.
 We give the case for two-sided jumps.  Note that this is not a restriction as
 the two-sided exit functionals  we consider are weakly continuous in the Skorokhod space. Therefore by taking limits as
 $\alpha(1-\rho)\rightarrow 1$ or $\alpha\rho\rightarrow 1$ we deduce
 identities for the case that $\xi^\uparrow$ is spectrally negative and
 spectrally positive respectively. Note that necessarily in the spectrally
 one sided case $\alpha\in(1,2)$.

\begin{theorem}Fix $\theta\geq 0$ and $-\infty<v<0<u<\infty$.
\begin{eqnarray*}
&&\hspace{-1cm}P\Big(\xi^\uparrow_{T^+_u}- u\in \ud\theta; T^+_u< T^-_v\Big) \\
&&=
\frac{\sin\pi\alpha(1-\rho)}{\pi} (e^u-1)^{\alpha(1-\rho)}(1-e^v)^{\alpha\rho}\\
&&\hspace{1cm}\times
 (e^{u+\theta} )^{\alpha\rho +1}(e^{u+\theta} -e^u)^{ -\alpha(1-\rho)}(e^{u+\theta}
 -e^v)^{-\alpha\rho}(e^{u+\theta}  -1)^{-1}\ud\theta
\end{eqnarray*}
and
\begin{eqnarray*}
&&\hspace{-1cm}P\Big(v-\xi^\uparrow_{T^-_v}\in \ud\theta; T^+_u> T^-_v\Big) \\
&&=\frac{\sin\pi\alpha\rho}{\pi} (1-e^v)^{\alpha\rho}(e^u -1)^{\alpha(1-\rho)}\\
&&\hspace{1cm}\times
 ( e^{v-\theta} )^{\alpha\rho+1}(e^{v} -e^{v-\theta})^{ -\alpha\rho}
 (e^u -e^{v-\theta})^{-\alpha(1-\rho)}(1-e^{v-\theta})^{-1}\ud\theta.
\end{eqnarray*}
\end{theorem}

\noindent{\it Proof.} Recall that $(X,\mathbf{P}_1)$ denotes an
$\alpha$-stable process initiated from $1$ and that
$(X,\mathbb{P}^\uparrow_1)$ is an $\alpha$-stable process conditioned to stay
positive initiated from $1$. From Lemma \ref{generic}, we have for
$\theta\geq 0$
\begin{eqnarray*}
&&\hspace{-1cm} P\Big(\xi^\uparrow_{T^+_u}\leq u+\theta; T^+_u< T^-_v\Big)\\
&&=\mathbb{P}^\uparrow_1\Big(X_{\sigma^+_{e^u}}\in[e^u, e^{u+\theta}]; \, \sigma^+_{e^u}< \sigma^-_{e^v}\Big)\\
&&=\int_{0}^{e^{u+\theta}
-e^u}(y+e^u)^{\alpha\rho}\mathbf{P}_1\Big(X_{\sigma^+_{e^u}} -e^u\in
\ud y;
\sigma^+_{e^u}< \sigma^-_{e^v}\Big)\\
&&=\int_{0}^{e^{u+\theta} -e^u}(y+e^u)^{\alpha\rho}\mathbf{P}_{1-e^v}
\Big(X_{\sigma^+_{(e^u-e^v)}} -(e^u-e^v)\in \ud y; \sigma^+_{(e^u-e^v)}< \sigma^-_{0}\Big)\\
&&=\frac{\sin\pi\alpha(1-\rho)}{\pi} (e^u-1)^{\alpha(1-\rho)}(1-e^v)^{\alpha\rho}\\
&&\hspace{1cm}\times\int_{0}^{e^{u+\theta} -e^u} (y+ e^u
)^{\alpha\rho}y^{ -\alpha(1-\rho)}(y+e^u -e^v)^{-\alpha\rho}(y+ e^u
-1)^{-1}\ud y
\end{eqnarray*}
from which the first part of the theorem follows.

The second part of the theorem can be proved in a similar way.
Indeed for $\theta \geq 0$
\begin{eqnarray*}
&&\hspace{-1cm} P\Big(\xi^\uparrow_{T^-_v}\geq v-\theta; T^+_u> T^-_v\Big)\\
&&=\mathbb{P}^\uparrow_1\Big(X_{\sigma^-_{e^v}}\in[e^{v-\theta}, e^{v}]; \, \sigma^+_{e^u}> \sigma^-_{e^v}\Big)\\
&&=\int_{0}^{e^{v}
-e^{v-\theta}}(e^v-y)^{\alpha\rho}\mathbf{P}_1\Big(e^v-X_{\sigma^-_{e^v}}
\in \ud y;
\sigma^+_{e^u}> \sigma^-_{e^v}\Big)\\
&&=\int_{0}^{e^{v} -e^{v-\theta}
}(e^v-y)^{\alpha\rho}\mathbf{P}_{1-e^v}\Big(-X_{\sigma^-_{0}} \in
\ud y;
\sigma^+_{(e^u-e^v)}> \sigma^-_{0}\Big)\\
&&=\frac{\sin\pi\alpha\rho}{\pi} (1-e^v)^{\alpha\rho}(e^u -1)^{\alpha(1-\rho)}\\
&&\hspace{1cm}\times\int_{0}^{e^{v} -e^{v-\theta}} ( e^v-y
)^{\alpha\rho}y^{ -\alpha\rho}(y+e^u
-e^v)^{-\alpha(1-\rho)}(y+1-e^v)^{-1}\ud y.
\end{eqnarray*}
This completes the proof.\QED

Note that since the process $(X,\mathbb{P}^\uparrow_x)$, $x>0$, is an $\alpha$-stable process
conditioned to stay positive, it follows that, in the case that
there are two-sided jumps, there is no creeping out of the interval
$(v,u)$ with probability one. That is to say.
\[
P\Big(\xi^\uparrow_{T^+_u}= u; T^+_u< T^-_v\Big)
=P\Big(\xi^\uparrow_{T^-_v}= v; T^+_u> T^-_v\Big) =0.
\]

Taking $v\downarrow-\infty$ in the first part of the above theorem
and $u\uparrow\infty$ in the second part we obtain the solution to
the one-sided exit problem as follows.
\begin{corollary}\label{onesideduparrow0}
Fix $\theta\geq 0$ and $-\infty<v<0<u<\infty$.
\begin{eqnarray*}
&&\hspace{-1cm}P\Big(\xi^\uparrow_{T^+_u}- u\in \ud\theta, T^+_u <\infty\Big) \\
&&= \frac{\sin\pi\alpha(1-\rho)}{\pi} (e^u-1)^{\alpha(1-\rho)}
 e^{u+\theta} (e^{u+\theta} -e^u)^{ -\alpha(1-\rho)}(e^{u+\theta}  -1)^{-1}\ud\theta
\end{eqnarray*}
and
\begin{eqnarray*}
&&\hspace{-1cm}P\Big(v-\xi^\uparrow_{T^-_v}\in \ud\theta;  T^-_v<\infty\Big) \\
&&=\frac{\sin\pi\alpha\rho}{\pi} (1-e^v)^{\alpha\rho}
 ( e^{v-\theta} )^{\alpha\rho+1}(e^{v}
 -e^{v-\theta})^{ -\alpha\rho}(1-e^{v-\theta})^{-1}\ud\theta.
\end{eqnarray*}
\end{corollary}

To give some credibility to these identities, and for future
reference, let us check that we may recover the identitiy in
Caballero and Chaumont \cite{CC} for the law of the minimum.

\begin{corollary} Let $\underline{\xi}^\uparrow_\infty=\displaystyle\inf_{t\ge
0}\xi^\uparrow_t$. For $z\geq 0$,
\[
P\Big(-\underline{\xi}^\uparrow_\infty \leq z\Big) =
(1-e^{-z})^{\alpha\rho}.
\]

\end{corollary}
\noindent{\it Proof}. The required probability may be identified as equal to $ P(T^-_{-z}=\infty)$ and hence, since there is no
probability of creeping over the level $-z$,
\begin{eqnarray*}
&&P\Big(-\underline{\xi}^\uparrow_\infty \leq z\Big)\\
&&= 1-\frac{\sin\pi\alpha\rho}{\pi} (1-e^{-z})^{\alpha\rho}
\int_0^\infty ( e^{-z-\theta} )^{\alpha\rho+1}(e^{-z}
-e^{-z-\theta})^{ -\alpha\rho}(1-e^{-z-\theta})^{-1}\ud\theta\\
&&=1-\frac{\sin\pi\alpha\rho}{\pi} (1-e^{-z})^{\alpha\rho}
\int_0^\infty ( e^{-\theta} )^{\alpha\rho+1}
(1-e^{-\theta})^{-\alpha\rho}(e^z-e^{-\theta})^{-1}\ud\theta.
\end{eqnarray*}
Next note that the integral in the right hand side satisfies
\begin{eqnarray}
&&\hspace{-1cm} \int_0^\infty ( e^{-\theta} )^{\alpha\rho+1}
(1 -e^{-\theta})^{ -\alpha\rho}(e^z-e^{-\theta})^{-1}\ud\theta\notag\\
&&=\int_1^\infty \frac{e^{-z}}{y(y-1)^{\alpha\rho} (y-e^{-z})}\ud y\notag\\
&&=\int_0^\infty \frac{e^{-z}}{(u+1)u^{\alpha\rho} (u+1-e^{-z})}\ud u\notag\\
&&=\int_0^\infty\left\{ \frac{1}{u^{\alpha\rho} (u+1-e^{-z})} -
\frac{1}{(u+1)u^{\alpha\rho} }\right\}\ud u\notag\\
&&=(1-e^{-z})^{-\alpha\rho}\int_0^\infty \frac{1}{v^{\alpha\rho}
(v+1)}dv - \int_0^\infty\frac{1}{(u+1)u^{\alpha\rho} }\ud u\notag\\
&&=[(1-e^{-z})^{-\alpha\rho} -1
]\int_0^\infty\frac{1}{(u+1)u^{\alpha\rho} }\ud u \label{ar-I}
\end{eqnarray}
where in the first equality we have applied the change of variable
$y=e^\theta$, in the second equality  $y=u+1$ and in the fourth
equality $u=(1-e^{-z})v$. Note also that by writing $w=(u+1)^{-1}$
we also discover that
\begin{equation}
 \int_0^\infty\frac{1}{(u+1)u^{\alpha\rho} }\ud u = \int_0^1
 (1-w)^{-\alpha\rho}w^{\alpha\rho-1}\ud w = \Gamma(1-\alpha\rho)\Gamma(\alpha\rho) =
 \frac{\pi}{\sin\pi\alpha\rho}.
 \label{ar-II}
\end{equation}
In conclusion we deduce that
\[
 \int_0^\infty ( e^{-\theta} )^{\alpha\rho+1}(1 -e^{-\theta})^{ -\alpha\rho}
 (e^z-e^{-\theta})^{-1}\ud\theta  =
 [(1-e^{-z})^{-\alpha\rho} -1 ] \frac{\pi}{\sin\pi\alpha\rho}
\]
and hence the required identity holds.\QED

Finally, to complete this subsection, when $(X,\p^\uparrow_1)$ is a spectrally
negative process we also gain some information concerning the scale
function, $W^{\uparrow, {\rm n}}$, of its underlying L\'evy process, denoted here by $\xi^{\uparrow,{\rm n}}$.
Specifically, in that case it is know that $1-\rho=1/\alpha$ (and
 $\alpha\in(1,2)$) and that $P(-\underline{\xi}^{\uparrow,{\rm n}}_\infty
\leq x) = mW^{\uparrow, {\rm n}}(x)$, where $m=E(\xi^{\uparrow,{\rm n}}_1)$. This implies
\[
 W^{\uparrow, {\rm n}}(x)(x) = \frac{1}{m}(1 - e^{-x})^{\alpha\rho}
 = \frac{1}{m}(1 - e^{-x})^{\alpha - 1}.
\]
 Recall that for a given spectrally negative L\'evy process
it is known that the Laplace transform of the scale function is
given by the inverse of the associated Laplace exponent (see for instance 
Theorem VII.8 in Bertoin \cite{Be}).   We
can therefore compute the Laplace exponent $\psi^\uparrow(\theta) =
\log E(e^{\theta \xi^{\uparrow,{\rm n}}_1})$ for $\theta\geq 0$, as follows:
\begin{eqnarray}
\psi^\uparrow(\theta) &=&m \left(\int_0^\infty e^{-\theta x}
(1- e^{-x})^{\alpha-1}\ud x\right)^{-1}\notag\\
&=&m\left(\int_0^1 u^{\theta - 1}(1-u)^{\alpha -1}\ud u\right)^{-1}\notag=m\frac{\Gamma(\theta+\alpha)}{\Gamma(\theta)\Gamma(\alpha)}.\label{-1}
\end{eqnarray}

The knowledge of the scale function allow us to write a stronger result than that given in
Corollary \ref{onesideduparrow0} as follows.

\begin{lemma}Let $\underline{\xi}^{\uparrow,{\rm n}}_t=\displaystyle\inf_{0\le s\le
t}\xi^{\uparrow,{\rm n}}_s$. For $v<0$,  $\theta\geq 0, \phi\geq \eta$ and
$\eta\in[0,-v]$ we have
 \begin{eqnarray*}
&& P\Big(v- \xi^{\uparrow,{\rm n}}_{T^-_v} \in \ud\theta, 
\xi^{\uparrow,{\rm n}}_{T^-_v - }-v\in
\ud\phi, \underline{ \xi}^{\uparrow,{\rm n}}_{T^-_v - }-v\in \ud\eta\Big)\\
 &&= K^{-1} \, (1 - e^{v+\eta})^{\alpha -
2}(e^{v+\eta})(e^{-\theta-\phi})^{\alpha}(1-e^{-\theta-\phi})^{-1-\alpha}\ud\theta
\ud\phi \ud\eta,
 \end{eqnarray*}
 where
 \[
 K=\frac{e^{(\alpha-2)v}}{\alpha(\alpha-1)}\int_1^{e^{-v}} \frac{(e^{-v}-y)}{y(y-1)^{\alpha-1}}\ud y-\frac{(1-e^v)^{\alpha-1}}{\alpha(\alpha-1)}\frac{\pi}{\sin \pi(\alpha-1)}.
 \]
\end{lemma}
\noindent{\it Proof}. First recall that the process $\xi^{\uparrow,{\rm n}}$
drifts towards $+\infty$ a.s. Taking this account, we have from
Example 8 of Doney and Kyprianou \cite{dk} that the required
probability is proportional to
\[
W^{\uparrow, {\rm n}}(-v-\ud\eta)\pi^\uparrow(-\theta -\phi)\ud\theta \ud\phi.
\]
Hence the triple law of interest has a density with respect to
$\ud\theta \ud\phi \ud\eta$ which is proportional to
\[ (1 -
e^{v+\eta})^{\alpha -
2}(e^{v+\eta})(e^{-\theta-\phi})^{\alpha}(1-e^{-\theta-\phi})^{-1-\alpha}.
\]
 For convenience let us write the constant of proportionality as $K^{-1}$. As
 $(X,\p^\uparrow_1)$ is derived from a spectrally negative stable process, it cannot
 creep downwards (cf. p175 of Bertoin \cite{Be}). This allows us to compute
 the unknown constant via the total probability formula and after a straightforward  computation, we have
\begin{eqnarray*}
K& =&\int_0^\infty\int_0^\infty\int_0^{-v}(1 -
e^{v+\eta})^{\alpha - 2}
(e^{v+\eta})(e^{-\theta-\phi})^{\alpha}(1-e^{-\theta-\phi})^{-1-\alpha}\ud\theta
\ud\phi \ud \eta\\
&=&\frac{e^{(\alpha-2)v}}{\alpha(\alpha-1)}\int_1^{e^{-v}} \frac{(e^{-v}-y)}{y(y-1)^{\alpha-1}}\ud y-\frac{(1-e^v)^{\alpha-1}}{\alpha(\alpha-1)}\frac{\pi}{\sin \pi(\alpha-1)}
\end{eqnarray*}
and the proof is complete.\QED


\subsection{Calculations for $\xi^*$}

Henceforth we shall assume that $(X,\mathbb{P}_x)$ is an $\alpha$-stable process
killed on first exit of the positive half line starting from $x>0$. As before, unless otherwise
stated, we shall assume that there are two-sided jumps, moreover,
spectrally one-sided results may be considered as limiting cases of
the two sided jumps case. We start with the two- and one-sided exit
problems with the latter as a limiting case of the former. We offer
no proof as the calculations are essentially the
same.

\begin{theorem}Fix $\theta\geq 0$ and $-\infty<v<0<u<\infty$.
\begin{eqnarray*}
&&\hspace{-1cm}P\Big(\xi^*_{T^+_u}- u\in \ud\theta; T^+_u< T^-_v\Big) \\
&&=
\frac{\sin\pi\alpha(1-\rho)}{\pi} (e^u-1)^{\alpha(1-\rho)}(1-e^v)^{\alpha\rho}\\
&&\hspace{1cm}\times
 (e^{u+\theta} )(e^{u+\theta} -e^u)^{ -\alpha(1-\rho)}(e^{u+\theta}
 -e^v)^{-\alpha\rho}(e^{u+\theta}  -1)^{-1}\ud\theta
\end{eqnarray*}
and
\begin{eqnarray*}
&&\hspace{-1cm}P\Big(v-\xi^*_{T^-_v}\in \ud\theta; T^+_u> T^-_v\Big) \\
&&=\frac{\sin\pi\alpha\rho}{\pi} (1-e^v)^{\alpha\rho}(e^u -1)^{\alpha(1-\rho)}\\
&&\hspace{1cm}\times
 ( e^{v-\theta} )(e^{v} -e^{v-\theta})^{ -\alpha\rho}
 (e^u -e^{v-\theta})^{-\alpha(1-\rho)}(1-e^{v-\theta})^{-1}\ud\theta.
\end{eqnarray*}
\end{theorem}

\begin{corollary}
Fix $\theta\geq 0$ and $-\infty<v<0<u<\infty$.
\begin{eqnarray*}
&&\hspace{-1cm}P\Big(\xi^*_{T^+_u}- u\in \ud\theta\Big) \\
&&= \frac{\sin\pi\alpha(1-\rho)}{\pi} (e^u-1)^{\alpha(1-\rho)}
( e^{u+\theta} )^{1-\alpha\rho}(e^{u+\theta} -e^u)^{ -\alpha(1-\rho)}(e^{u+\theta}  -1)^{-1}\ud\theta
\end{eqnarray*}
and
\begin{eqnarray*}
&&\hspace{-1cm}P\Big(v-\xi^*_{T^-_v}\in \ud\theta;  T^-_v<\infty\Big) \\
&&=\frac{\sin\pi\alpha\rho}{\pi} (1-e^v)^{\alpha\rho}
 ( e^{v-\theta} )(e^{v}
 -e^{v-\theta})^{ -\alpha\rho}(1-e^{v-\theta})^{-1}\ud\theta.
\end{eqnarray*}
\end{corollary}

One may think of computing the distribution of the maximum, $\overline{\xi}^*_\infty$, and the minimum,
$\underline{\xi}^*_\infty$, of $\xi^*$ in a similar way to the previous section by integrating out $u$
and $v$ in the above corollary. The law of the minimum was already computed in Caballero and Chaumont
(2006) and we refrain from producing the alternative computations here.
For the maximum, an easier approach is at hand. Since $\xi^*$ is derived from a stable process killed
on first exit of the positive half line one may write
\[
P\Big(\overline{\xi}^*_\infty \leq z\Big)=P\Big(\exp\{\overline{\xi}^*_\infty\} \leq e^z\Big) =
\mathbf{P}_1(\sigma^+_{e^z} > \sigma^-_0)  = \mathbf{P}_{e^{-z}}(\sigma^+_1>\sigma^-_0).
\]
The probability on the right hand side above may be obtained from Theorem \ref{rogozin} by a
straightforward integration. The latter calculation has already been performed however in Rogozin
\cite{ro} and is equal to
\[
\frac{\Gamma(\alpha)}{ \Gamma(\alpha\rho)\Gamma(\alpha(1-\rho))}\int_0^{1-e^{-z}} y^{\alpha\rho-1}
(1-y)^{\alpha (1-\rho)  -1}\ud y
\]
Hence, together with the result for the minimum from Caballero and Chaumont (2006) which we include
for completeness, we have the following corollary.

\begin{corollary}
For $z \geq 0$ we have that
\[
P\Big(\overline{\xi}^*_\infty \in dz\Big)= \frac{\Gamma(\alpha)}{\Gamma(\alpha\rho)\Gamma(\alpha(1-\rho))}
(e^{-z})^{\alpha(1-\rho)}(1- e^{-z})^{\alpha\rho - 1}dz
\]
and
\[
P\Big(-\underline{\xi}^*_\infty \in dz\Big) = \frac{1}{\Gamma(\alpha\rho)\Gamma(1-\alpha\rho)}
(e^z-1)^{\alpha\rho}dz.
\]
\end{corollary}

In the case that $\xi^*$ is spectrally one sided, it seems difficult
to use the above result to extract information about any underlying
scale functions. The reason for this is that the process $\xi^*$ is
exponentially killed at a rate which is intimately linked to its
underlying parameters and not at a rate which can be independently
varied.

\subsection{Calculations for $\xi^\downarrow$}

Henceforth we shall assume that $(X,\mathbb{P}^\downarrow_x)$ is an $\alpha$-stable process
conditioned to hit zero continuously starting from $x>0$. Again, unless otherwise
stated, we shall assume that there are two-sided jumps, and
spectrally one-sided results may be considered as limiting cases of
the two sided jumps case. We follow the same programme as the
previous two sections dealing with the two- and one-sided exit
problems without offering proofs since they follow from the calculations for $\xi^\uparrow$ and Proposition \ref{4590}

\begin{theorem}Fix $\theta \geq 0$.
\begin{eqnarray*}
&&\hspace{-1cm}P\Big(\xi^\downarrow_{T^+_u}- u\in \ud\theta; T^+_u< T^-_v\Big) \\
&&=  \frac{\sin\pi\alpha(1-\rho)}{\pi} (e^u-1)^{\alpha(1-\rho)}(1-e^v)^{\alpha\rho}\\
&&\hspace{1cm}\times
 (e^{u+\theta} )^{\alpha\rho}(e^{u+\theta} -e^u)^{ -\alpha(1-\rho)}
 (e^{u+\theta}  -e^v)^{-\alpha\rho}(e^{u+\theta}  -1)^{-1}\ud\theta
\end{eqnarray*}
and
\begin{eqnarray*}
&&\hspace{-1cm}P\Big(v-\xi^\downarrow_{T^-_v}\in \ud\theta; T^+_u> T^-_v\Big) \\
&&=\frac{\sin\pi\alpha(1-\rho)}{\pi} (e^u-1)^{\alpha(1-\rho)}(1-e^v)^{\alpha\rho}\\
&&\hspace{1cm}\times
 ( e^{v-\theta} )^{\alpha\rho}(e^{v} -e^{v-\theta})^{ -\alpha\rho}
 (e^u -e^{v-\theta})^{-\alpha(1-\rho)}(1-e^{v-\theta})^{-1}\ud\theta.
\end{eqnarray*}
\end{theorem}

\begin{corollary}\label{onesideduparrow}
Fix $\theta\geq 0$.
\begin{eqnarray*}
&&\hspace{-1cm}P\Big(\xi^\uparrow_{T^+_u}- u\in \ud\theta; T^+_u<\infty\Big) \\
&&= \frac{\sin\pi\alpha(1-\rho)}{\pi} (e^u-1)^{\alpha(1-\rho)}
(e^{u+\theta} -e^u)^{ -\alpha(1-\rho)}(e^{u+\theta}  -1)^{-1}\ud\theta
\end{eqnarray*}
and
\begin{eqnarray*}
&&\hspace{-1cm}P\Big(v-\xi^\uparrow_{T^-_v}\in \ud\theta;  T^-_v<\infty\Big) \\
&&=\frac{\sin\pi\alpha\rho}{\pi} (1-e^v)^{\alpha\rho}
 ( e^{v-\theta} )^{\alpha\rho}(e^{v} -e^{v-\theta})^{ -\alpha\rho}(1-e^{v-\theta})^{-1}\ud\theta.
\end{eqnarray*}
\end{corollary}

From the above corollary we proceed to obtain the law of the maximum
of $\xi^\downarrow$ (recalling that it is a process with drift to $-\infty$).

\begin{corollary}
For $z\geq 0$
\[
P\Big(\overline{\xi}^\downarrow_\infty \leq z\Big) =
(1-e^{-z})^{\alpha(1-\rho)}
\]
\end{corollary}
\noindent{\it Proof}. Similarly to the calculations in  Section \ref{xiup} we make use of
the fact that 
\[
P\Big(\overline{\xi}^\downarrow_\infty \leq z\Big) =
P(T^+_z =\infty).
\]
Hence
\begin{eqnarray*}
&&\hspace{-1cm}P\Big(\overline{\xi}^\downarrow_\infty \leq z\Big) \\
&&=1 - \frac{\sin\pi\alpha(1-\rho)}{\pi}
(e^z-1)^{\alpha(1-\rho)}\int_0^\infty
(e^{z+\theta} -e^z)^{ -\alpha(1-\rho)}(e^{z+\theta}  -1)^{-1}\ud\theta\\
&&1 - \frac{\sin\pi\alpha(1-\rho)}{\pi} (1-e^{-z})^{\alpha(1-\rho)}
\int_0^\infty (e^{-\theta})^{\alpha(1-\rho)+1} (1-e^{-\theta} )^{
-\alpha(1-\rho)}(e^{z}  -e^{-\theta})^{-1}\ud\theta.
\end{eqnarray*}
Next note that the integral on the right hand side above has been
seen before in (\ref{ar-I}) except for the case that $\rho$ is
replaced by $1-\rho$. We thus obtain from (\ref{ar-I}) and
(\ref{ar-II})
\[
\int_0^\infty (e^{-\theta})^{\alpha(1-\rho)+1} (1-e^{-\theta} )^{
-\alpha(1-\rho)}(e^{z}  -e^{-\theta})^{-1}d\theta
=[(1-e^{-z})^{-\alpha(1-\rho)} -1]\frac{\pi}{\sin\pi\alpha(1-\rho)}
\]
and hence
\[
P\Big(\overline{\xi}^\downarrow_\infty \leq z\Big)
=(1-e^{-z})^{\alpha(1-\rho)}
\]
as required.\QED

Now, we suppose  that $(X,\p^{\downarrow}_1)$ has only positve jumps, in which case  $\rho  =
1/\alpha$. We denote by $\xi^{\downarrow,{\rm p}}$ its underlying L\'evy process in this particular case. The associated scale function of $\xi^{\downarrow,{\rm p}}$ can be identified as
\begin{equation}\label{duality}
E(-\xi^{\downarrow,{\rm p}}_1)W^{\downarrow,{\rm p}}(x) = P\Big(\overline{\xi}^{\downarrow,{\rm p}}_\infty \leq x\Big) = (1-
e^{-x})^{\alpha-1}=W^{\uparrow, {\rm n}}(x)E(\xi^{\uparrow, {\rm n}}_1),
\end{equation}
where $W^{\uparrow, {\rm n}}$ is the scale function of the spectrally negative L\' evy process 
$\xi^{\uparrow,{\rm n}}$. \\
This observation reflects the duality property for positive
self-similar Markov processes in this particular case (see section 2 in Bertoin and Yor \cite{BeY}).
More precisely, we have the  duality property between the resolvent operators of $(X,\p^{\uparrow}_{x})$, when it is spectrally negative, and 
 $(X,\p^{\downarrow}_{x})$, for $x>0$.

From the identification of the scale function in (\ref{duality}) and Lemma 2, it is possible to give a triple law for the process $\xi^{\downarrow,{\rm p}}$ at first passage over the level $x>0$. 
\begin{lemma}Let $\overline{\xi}^{\downarrow,{\rm p}}_t=\displaystyle\sup_{0\le s\le
t}\xi^{\downarrow,{\rm p}}_s$. For $x>0$,  $\theta\geq 0, \phi\geq \eta$ and
$\eta\in[0,x]$ we have
 \begin{eqnarray*}
&& P\Big(\xi^{\downarrow,{\rm p}}_{T^+_x}-x \in \ud\theta, 
x-\xi^{\downarrow,{\rm p}}_{T^+_x - }\in
\ud\phi, x-\overline{ \xi}^{\downarrow, {\rm p}}_{T^+_x - }\in \ud\eta\Big)\\
 &&= K^{-1} \, (1 - e^{-x+\eta})^{\alpha -
2}e^{-x+\eta}e^{\theta+\phi}(e^{\theta+\phi}-1)^{-1-\alpha}\ud\theta
\ud\phi \ud\eta,
 \end{eqnarray*}
 where
 \[
 K=\frac{e^{(\alpha-2)v}}{\alpha(\alpha-1)}\int_1^{e^{-v}} \frac{(e^{-v}-y)}{y(y-1)^{\alpha-1}}\ud y-\frac{(1-e^v)^{\alpha-1}}{\alpha(\alpha-1)}\frac{\pi}{\sin \pi(\alpha-1)}.
 \]
\end{lemma}

In the remainder of this section, we assume that $(X,\p^{\downarrow}_x)$ is spectrally negative and we denote its underlying L\'evy process by $\xi^{\downarrow, {\rm n}}$ for its underlying L\'evy processes.  The identification of the scale function of the L\'evy process
$\xi^{\uparrow,{\rm n}}$ and Proposition \ref{4590} inspires the following result
which identifies the scale function of the L\'evy process
$\xi^{\downarrow, {\rm n}}$. We emphasized that $\xi^{\downarrow,{\rm n}}$ is in fact
$\xi^{\uparrow,{\rm n}}$ conditioned to drift to $-\infty$.

\begin{lemma} The Laplace exponent of 
$\xi^{\downarrow,{\rm n}}$ satisfies
\[
\psi^{\downarrow}(\theta)
=m\frac{\Gamma(\theta-1+\alpha)}{\Gamma(\theta-1)\Gamma(\alpha)}
\]
for $\theta\geq 0$, and where $m=E(-\xi^{\downarrow,{\rm n}})$. Moreover, its associated scale function may be identified as
\[
W^{\downarrow, {\rm n}}(x) = \frac{1}{m}(1-e^{-x})^{\alpha
-1}e^{ x}.
\]
\end{lemma}
\noindent{\it Proof}.   From the two-sided exit problem for spectrally negative L\'evy processes (see for instance Chapter VII in Bertoin \cite{Be}), we know that 
\[
P\Big(\underline{\xi}^{\uparrow, {\rm n}}_{T^{+}_y}>-x\Big)=\frac{W^{\uparrow, {\rm n}}(x)}{W^{\uparrow, {\rm n}}(x+y)}\qquad\textrm{ for }\,\, x,y>0.
\]
On the other hand, from Proposition 1, we get that
\[
P\Big(\underline{\xi}^{\uparrow, {\rm n}}_{T^{+}_y}>-x\Big)=e^yP\Big(\underline{\xi}^{\downarrow, {\rm n}}_{T^{+}_y}>-x\Big)=e^y\frac{W^{\downarrow, {\rm n}}(x)}{W^{\downarrow,{\rm n}}(x+y)}.
\]
Hence, from the form of the scale function of $\xi^{\uparrow,{\rm n}}$,  it follows that
\[
P\Big(\underline{\xi}^{\downarrow, {\rm n}}_{T^{+}_y}>-x\Big)=\frac{W^{\downarrow, {\rm n}}(x)}{W^{\downarrow,{\rm n}}(x+y)}=\frac{e^x(1-e^{-x})^{\alpha-1}}{e^{x+y}(1-e^{-x-y})^{\alpha-1}}.
\]
Taking $y$ to $\infty$, one deduces
\[
W^{\downarrow, {\rm n}}(x)=\frac{1}{m}e^x(1-e^{-x})^{\alpha-1}.
\]
Finally, since
\[
\int_0^\infty e^{-\theta x}e^{ x}(1-e^{-x})^{\alpha-1}\ud x =
\int_{0}^1 u^{(\theta-1)-1}(1-u)^{\alpha-1}\ud u =
\frac{\Gamma(\theta-1)\Gamma(\alpha)}{\Gamma(\theta-1+\alpha)},
\]
for $\theta>1$, it is clear from the Laplace transform (cf Chapter 8 of
Kyprianou \cite{Ky}) of the scale function that 
\[
\psi^{\downarrow} (\theta)=m\frac{\Gamma(\theta-1+\alpha)}{\Gamma(\theta-1)\Gamma(\alpha)},\, \theta\geq 0,
\] 
is the associated
Laplace exponent of $\xi^{\downarrow, {\rm n}}$.\QED

One may also write down a triple law for the first passage problem of $\xi^{\downarrow,{\rm n}}$ 
as we have seen before for $\xi^{\uparrow,{\rm n}}$.

\begin{lemma}Let $\underline{\xi}^{\downarrow,{\rm n}}_t=\displaystyle\inf_{0\le s\le
t}\xi^{\downarrow,{\rm n}}_s$. For $v<0$,  $\theta\geq 0, \phi\geq \eta$ and
$\eta\in[0,-v]$ we have
 \begin{eqnarray*}
&& P\Big(v- \xi^{\downarrow,{\rm n}}_{T^-_v} \in \ud\theta, 
\xi^{\downarrow,{\rm n}}_{T^-_v - }-v\in
\ud\phi, \underline{ \xi}^{\downarrow,{\rm n}}_{T^-_v - }-v\in \ud\eta\Big)\\
 &&= K^{-1} \, e^{-q(\phi-\eta)}(e^{-(v+\eta)}+\alpha-2)(1 - e^{v+\eta})^{\alpha -
2}(e^{-\theta-\phi})^{\alpha}(1-e^{-\theta-\phi})^{-1-\alpha}\ud\theta
\ud\phi \ud\eta,
 \end{eqnarray*}
 where $q>0$ and 
 \[
 K=\frac{1}{\alpha}\int_0^{-v} \int_\eta^{\infty} e^{-q(\eta-\phi)}(e^{-v-\eta}+\alpha-2)(1-e^{v+\eta})^{\alpha-2}(e^{\phi}-1)^{-\alpha}\ud \phi\ud \eta.
 \]
\end{lemma}
\noindent{\it Proof}. First recall that the process $\xi^{\downarrow,{\rm n}}$
drifts towards $-\infty$ a.s. Again, we have from
Example 8 of Doney and Kyprianou \cite{dk} that the required
probability is proportional to
\[
e^{-q(\phi-\eta)}W^{\downarrow, {\rm n}}(-v-\ud\eta)\pi^{\downarrow}(-\theta -\phi)\ud\theta \ud\phi,
\]
where $q>0$ is the killing rate of the descending ladder height process (see for instance Chapter VI in Bertoin \cite{Be} for a proper definition) of $\xi^{\downarrow,{\rm n}}$.\\ 
Hence the triple law of interest has a density with respect to
$\ud\theta \ud\phi \ud\eta$ which is proportional to
\[ e^{-q(\phi-\eta)}(1 -
e^{v+\eta})^{\alpha -
2}(e^{-(v+\eta)} +\alpha-2)(e^{-\theta-\phi})^{\alpha}(1-e^{-\theta-\phi})^{-1-\alpha}.
\]
 As $(X,\p^\downarrow_1)$ is derived from a spectrally negative stable process, it cannot
 creep downwards (cf. p175 of Bertoin \cite{Be}). This allows us to compute
 the unknown constant of proportionality $K^{-1}$ via the total probability formula and after a straightforward  computation, we have
\begin{eqnarray*}
K& =&\int_0^\infty\int_0^\infty\int_0^{-v}e^{-q(\phi-\eta)}(1 -
e^{v+\eta})^{\alpha - 2}
(e^{v+\eta})(e^{-\theta-\phi})^{\alpha}(1-e^{-\theta-\phi})^{-1-\alpha}\ud\theta
\ud\phi \ud \eta\\
&=&\frac{1}{\alpha}\int_0^{-v} \int_\eta^{\infty} e^{-q(\eta-\phi)}(e^{-v-\eta}+\alpha-2)(1-e^{v+\eta})^{\alpha-2}(e^{\phi}-1)^{-\alpha}\ud \phi\ud \eta.
\end{eqnarray*}
and the proof is complete.\QED


\section{Entrance laws for L\'evy-Lamperti processes: points}

In this section we explore the two-point hitting problem for the
L\'evy-Lamperti processes $\xi^\uparrow$  and
$\xi^{\downarrow}$. There has been little work dedicated to this
theme in the past with the  paper of Getoor \cite{Get} being our
principle reference.

Henceforth we shall denote by  $(X,\mathbf{P}_x)$ a {\it symmetric}
$\alpha$-stable process issued from $x>0$ where $\alpha\in(1,2)$. An
important quantity in the forthcoming analysis is the resolvent density  of
the process $(X,\mathbf{P}_x)$ killed on exiting $(0,\infty)$. The latter is known
to have a density
\[
 u(x,y)\ud y = \int_0^\infty \ud t \cdot \mathbf{P}_x(X_t \in dy, t< \sigma^-_0)
\]
for $x,y>0$. From Blumenthal et al. \cite{Bl} we know that
\[
 \int_0^\infty \ud t\cdot \mathbb{P}_x(X_t \in \ud y, t< \sigma^+_a\wedge \sigma^-_0)
= \left\{ \frac{|x- y|^{\alpha
-1}}{2^\alpha \Gamma(\alpha/2)} \int_0^{s(x,y,a)} \frac{u^{\alpha/2 -1}}{(u+1)^{1/2}}\ud u
\right\}\ud y
\]
where
\[
 s(x,y,a) = \frac{4xy}{(x-y)^2}\frac{(a-x)(a-y)}{a^2}.
\]
It now follows taking limits as $a\uparrow\infty$
that
\[
u(x,y)= \frac{1}{2^\alpha \Gamma(\alpha/2)} |x- y|^{\alpha
-1}\int_0^{4xy/(x-y)^2} \frac{u^{\alpha/2 -1}}{(u+1)^{1/2}}\ud u.
\]
According to the method presented in Getoor \cite{Get} one may compute 
\[
\mathbf{P}_x(X_{\sigma_{\{a,b\}}} =a; \, \sigma_{\{a,b\}} < \sigma^-_0) 
\]
where $\sigma_{\{a,b\}}= \inf\{t>0 : X_t = a \text{ or }b\}$ and $a,b>0$
using the following technique. The two point hitting probability in Getoor \cite{Get} is
given by the formula
\[
\mathbf{P}_x(X_{\sigma_{\{a,b\}}} =a; \, \sigma_{\{a,b\}} < \sigma^-_0) = -
\frac{q(x,a)}{q(x,x)} \label{hittingprob}
\]
where the $\{x,a,b\}\times\{x,a,b\}$-matrix $Q$ is defined by
\[
Q= - U^{-1} \label{inverse}
\]
and the $\{x,a,b\}\times\{x,a,b\}$-matrix is given by
\[
U= \left(
\begin{array}{ccc}
u(x,x) & u(x,a) & u(x,b)\\
u(a,x) & u(a,a) & u(a,b) \\
u(b,x) & u(b,a) & u(b,b)
 \end{array}
\right).
\]
In particular an easy computation shows that
\begin{equation}
 \mathbf{P}_x(X_{\sigma_{\{a,b\}}} =a;  \, \sigma_{\{a,b\}} < \sigma^-_0) =\frac{\frac{u(x,a)}{u(b,a)} - \frac{u(x,b)}{u(b,b)}}{\frac{u(a,a)}{u(b,a)} - \frac{u(a,b)}{u(b,b)}}
\label{hittingprob}
\end{equation}

Recalling the definitions of $\xi^\uparrow$ and $\xi^\downarrow$ as
the L\'evy-Lamperti processes associated now with our symmetric
stable process conditioned to stay positive and conditioned to be
killed continuously at the origin respectively we obtain the
following result.
\begin{theorem}
Fix $\alpha\in(1,2)$ and $-\infty<v<0<u<\infty$. Define
\[
T_{\{v,u\}}=\inf\{t>0 : \xi_t\in\{v,u\}\}
\]
where $\xi$ plays the role of either $\xi^\uparrow$ or $\xi^\downarrow$.
 We have
\[
 P\Big(\xi^\uparrow_{   T_{\{v,u\}}  } =v\Big) = (e^v)^{\alpha/2} f(1, e^v, e^u)
\]
and
\[
 P\Big(\xi^\downarrow_{T_{\{v,u\}}} =v\Big) = (e^v)^{\alpha/2-1} f(1, e^v, e^u)
 \]
where
\[
f(x,a,b) = \frac{\frac{u(x,a)}{u(b,a)} - \frac{u(x,b)}{u(b,b)}}{\frac{u(a,a)}{u(b,a)} - \frac{u(a,b)}{u(b,b)}}.
\]
\end{theorem}


\section{Exponential functionals of L\'evy-Lamperti processes.}

We begin this section by recalling a crucial expression for the
entrance law at 0 of pssMp's. In \cite{beC, BeY}, the authors proved
that if a non arithmetic L\'evy process $\xi$ satisfies
$E(|\xi_1|)<\infty$ and $0<E(\xi_1)<+\infty$, then its corresponding
pssMp $(X,\p_{x})$ in the Lamperti representation converges weakly
as $x$ tends to $0$, in the sense of finite dimensional
distributions towards a non degenerated probability law $\p_0$.
Under these conditions, the entrance law under $\p_0$ is described
as follows: for every $t>0$ and every measurable function $f:\R_+\to
\R_+$,
\begin{equation}\label{entlaw}
\e_0\big(f(X_t)\big)=\frac{1}{\alpha E(\xi_1)}E\left(I(\xi)^{-1}
f\big(tI(\xi)^{-1}\big)\right)\,,
\end{equation}
where $I(\xi)$ is the exponential functional:
\[I(\xi)=\int_{0}^{\infty} \exp\{-\alpha \xi_{s}\}\,
\ud s.\] Necessary and sufficient conditions for the weak
convergence of $(X, \p_x)$ on the Skorokhod's space were given in
\cite{CCh}. Recall that $(X, \px^{\uparrow})$ denotes a stable
L\'evy process conditioned to stay positive as it has been defined
in section \ref{prelim}. Then we easily check that $\xi^\uparrow$
satisfies conditions for the weak convergence of $(X,\p_x^{\uparrow})$ given in
\cite{beC, BeY, CCh}. Note also that in this particular case, the
weak convergence of $(X,\p^{\uparrow}_x)$ had been proved in a direct way in
\cite{Ch}. We denote the limit law by $\p^{\uparrow}$.\\

\noindent We first investigate the tail behaviour of the law of
$I(\xi^\uparrow)$.
\begin{theorem}\label{asymptotic}
The law of  $I(\xi^\uparrow)$  is absolutely continuous with respect
to the Lebesgue measure. The density of $I(\xi^\uparrow)^{-1}$ is given by:
\begin{equation}\label{elaw}
P\Big(I(\xi^\uparrow)^{-1}\in \ud y\Big)=\alpha E(\xi^{\uparrow}_1)y^{\alpha\rho -1}q_1(y)\ud y,
\end{equation}
where $q_t$ is the density of  the entrance law of the excursion measure of the reflected process $(X-\underline{X}, \mathbf{P}_0)$, where $\underline{X}_t=\inf_{0\leq s\leq t}X_s$. Moreover, the law of $I(\xi^\uparrow)$ behaves as
\begin{equation} \label{at}
P(I(\xi^{\uparrow})\ge x)\sim C_1
x^{-\alpha}\,,\;\;\mbox{as $x\rightarrow+\infty$.}
\end{equation}
If $X$ has  positive jumps, then
\begin{equation}\label{at0}
P(I(\xi^{\uparrow})\le x)\sim C_2 x^{\alpha(\rho-1)-1}\,,\;\mbox{as
$x\rightarrow0$.}
\end{equation}
The constants $C_1$ and $C_2$  depend only on $\alpha$ and $\rho$.
\end{theorem}
\noindent In the case where the process has positive jumps, the  law
of $I(\xi^\uparrow)$ is given explicitly in the next theorem.\\

\noindent {\it Proof}. Let $n$ be the measure of the excursions away
form 0 of the reflected process $X-\underline{X}$ under
$\mathbf{P}_0$. It is proved in \cite{MS} that the entrance law of
$n$ is absolutely continuous with respect to the Lebesgue measure.
Let us denote by $q_t$ its density. Then from \cite{Ch}, the entrance law of 
$(X,\p^{\uparrow})$ is
related to $q_1$ by
\[\p^{\uparrow}(X_1\in \ud y)=y^{\alpha\rho}q_1(y)\ud y\,,\;\;\;y\ge0\,.\]
We readily derive (\ref{elaw}) from  identity
(\ref{entlaw}). Moreover from  (3.18) in
\cite{MS}:
\[\int_0^x q_1(y)\,\ud y\sim Cx^{\alpha(1-\rho)+1}\,,\;\;\;\mbox{as $x\rightarrow0$,}\]
and from  (3.20) of the same paper, if $X$ has positive
jumps, then:
\[\int_x^\infty q_1(y)\,\ud y\sim C'x^{-\alpha}\,,\;\;\;\mbox{as $x\rightarrow+\infty$.}\]
This together with (\ref{entlaw}) implies (\ref{at}) and (\ref{at0}). The constants $C$ and 
$C'$ depend only on $\alpha$ and $\rho$. \QED
\noindent Another way to prove part $(i)$ of this theorem is to use
a result due to M\'ejane \cite{Me} and Rivero \cite{ri1} which
asserts that for a non arithmetic L\'evy process $\xi$, if Cramer's
condition is satisfied for $\theta>0$, i.e. $E(\exp\theta\xi_1)=1$
and $E(\xi_1^+\exp{\theta\xi_1})<\infty$, then $P(I(\xi)\ge x)\sim
Cx^{-\alpha\theta}$. These arguments and Proposition \ref{4590}
allow us to obtain the asymptotic behaviour at $+\infty$ of
$P(I(-\xi^{\downarrow})\ge x)$:

\begin{proposition}\label{downasymptotic}
The law of  $I(-\xi^\downarrow)$ behaves as
\begin{equation}
P(I(-\xi^{\downarrow})\ge x)\sim C_3 x^{-\alpha}\,,\;\;\mbox{as
$x\rightarrow+\infty$.}
\end{equation}
The constant $C_3$ depends only on $\alpha$ and $\rho$.
\end{proposition}

\noindent Now we consider the exponential functional
\[I(-\xi^*)=\int_{0}^{\infty} \exp\{\alpha\xi^{*}_s\}\, \ud s.\]
Recall from section \ref{prelim} that $\xi^*$ is the L\'evy process
which is associated to the pssMp $(X_t\ind_{\{t<T\}},\mathbf{P}_x)$
by the Lamperti representation. From this representation, we may
check path by path the equality 
\[x^\alpha I(-\xi^*)=T\,.\]
Moreover, it follows from Lemma 1 in \cite{CC} that when
$(X,\mathbf{P}_x)$ has negative jumps, $\mathbf{P}_x(T\le t)\sim
\frac{c_-t}{\alpha x^\alpha}$, as $t$ tends to 0. This result  leads
to:
\begin{proposition}\label{42361}
Suppose that $\xi^*$ has negative jumps, then the law of $I(-\xi^*)$
behaves as
\begin{equation}
P(I(-\xi^*)\le x)\sim \frac{c_-}{\alpha} x^{-1}\,,\;\;\mbox{as
$x\rightarrow0$.}
\end{equation}
\end{proposition}

\noindent In the remainder of this section, we assume that $(X, \px^{\uparrow})$ has no positive jumps.

\begin{theorem} The law of exponential functional $I(\xi^{\uparrow,{\rm n}})=
\int_{0}^{\infty} \exp\{-\alpha\xi^{\uparrow,{\rm n}}_{s}\}\, \ud s$ is
absolutely continuous with respect to the Lebesgue measure and has a
continuous density $p^{\uparrow,{\rm n}}(\cdot)$ which has the
following representation by power series
\[
p^{\uparrow, {\rm n}}(x)=-\frac{c^{-1}}{\pi
x}\sum_{n=1}^{\infty}\Gamma\left(1+\frac{n}{\alpha}\right)\sin\left(\frac{\pi
n}{\alpha}\right)\frac{(-x^{-1/\alpha})^n}{n!}, \qquad \textrm{for
}\quad x>0,
\]
where $c=c_{-}\Gamma(2-\alpha)\alpha^{-1}(\alpha-1)^{-1}>0$.\\
Moreover the positive entire moments of $(X,\p^{\uparrow})$, for $t>0$, are
given by the identity
\begin{equation}\label{moments}
\e^{\uparrow}\Big(\big(X_t\big)^k\Big)=(mt)^k
\frac{\Gamma(\alpha(k+1))}{\Gamma(\alpha)^{(k+1)}k!}, \qquad k\geq
1,
\end{equation}
and its law is absolutely continuous with respect to the Lebesgue
measure and has a continuous density $p_t(\cdot)$ which has the
following representation by power series
\[
p_t(x)=-\frac{t^{-1/\alpha}(\alpha-1)}{\Gamma(2-\alpha)c_{-} m \pi}
\sum_{n=1}^{\infty}\Gamma\left(1+\frac{n}{\alpha}\right)\sin\left(\frac{\pi
n}{\alpha}\right)\frac{(-x^{1/\alpha})^n}{n!}, \qquad\textrm{for
}\quad  x>0.
\]
\end{theorem}
\noindent {\it Proof}. From  Bertoin and Yor \cite{BeY2} , we know
that the distribution of the exponential functional of a spectrally negative L\'evy
process is determined by its negative entire
moments. In particular, the exponential functional for
$\xi^{\uparrow,{\rm n}}$  satisfies
\begin{equation}\label{nmom}
E\left(I(\xi^{\uparrow, {\rm n}})^{-k}\right)= \alpha
m^k\frac{\Gamma(k\alpha)}{\Gamma(\alpha)^{k}(k-1)!},
\end{equation}
with the convention that the right-hand side equals $\alpha m$ for
$k=1$. In particular, from (\ref{entlaw}) we have that
\[
\e^{\uparrow}\Big(\big(X_t\big)^k\Big)=(mt)^k
\frac{\Gamma(\alpha(k+1))}{\Gamma(\alpha)^{(k+1)}k!},
\]
which proves the identity (\ref{moments}).\\
Now,  from the time reversal property of Theorem VII.18 in
\cite{Be}, we deduce that the last passage time of $(X,\p^{\uparrow})$,
defined by
 \[
 U_{x}=\sup\Big\{t\geq 0:\, X_t\leq x\Big\}\qquad \textrm{for }\quad x\geq 0,
 \]
 is a stable subordinator of index $1/\alpha$. More precisely, its Laplace exponent
 is given by $\Phi(\lambda)=(\lambda/c)^{1/\alpha}$, where
 $c=c_{-}\Gamma(2-\alpha)\alpha^{-1}(\alpha-1)^{-1}$. According to Zolotarev
 \cite{Zo},
 stable subordinators have  continuous densities with respect  to the Lebesgue
 measure which  may be represented by power series. More precisely,
the density of a normalized stable subordinator of index
$\beta\in(0,1)$, i.e. $\Phi(\lambda)=\lambda^{\beta}$,  is given by
\begin{equation}\label{subd}
\rho_t(x, \beta)=-\frac{t^{-1/\beta}}{\pi
x}\sum_{n=1}^{\infty}\Gamma(1+n\beta)\sin(\beta\pi
n)\frac{(-x^{-\beta})^n}{n!}, \qquad \textrm{for }\quad x>0.
\end{equation}
On the other hand  from Proposition 1 in \cite{CP}, we know that
$U_x$ has the same law as $x^{\alpha}I(\xi^{\uparrow, {\rm n}})$. Hence,
$I(\xi^{\uparrow, {\rm n}})$ satisfies that
\[
E\Big(\exp\big\{-\lambda
I(\xi^{\uparrow, {\rm n}})\big\}\Big)=e^{-(\lambda/c)^{1/\alpha}}, \qquad
\lambda\geq 0,
\]
and its density  $p^{\uparrow, {\rm n}}(x)$ is given by
$ \rho_{c^{1/\alpha}}(x,1/\alpha)$, which proves the first part of the theorem.\\
Finally from (\ref{entlaw}), we deduce that the density of
$X^{(0)}_1$ is given by
\[
p_1(x)=-\frac{c^{-1}}{\alpha m \pi}
\sum_{n=1}^{\infty}\Gamma\left(1+\frac{n}{\alpha}\right)\sin\left(\frac{\pi
n}{\alpha}\right)\frac{(-x^{1/\alpha})^n}{n!}, \qquad\textrm{for
}\quad  x>0.
\]
The proof is now complete.\QED

\begin{theorem} The exponential functional
$I(-\xi^{*})=\int_{0}^{\infty} \exp\{\alpha\xi^{*}_s\}\, \ud s$ has
a  continuous density $p^*(\cdot)$ with respect to the Lebesgue
measure which has the following representation by power series
\[
p^{*}(x)=c^{1/\alpha}\sum_{n=1}^{\infty}\frac{\alpha
n-1}{\Gamma(\alpha
n)\Gamma(-n+1+1/\alpha)}x^{\alpha(2-n\alpha)},\quad
\textrm{for}\quad x>0,\] where
$c=c_{+}\Gamma(2-\alpha)\alpha^{-1}(\alpha-1)^{-1}>0.$
\end{theorem}
\noindent {\it Proof:} First, let us  define $\hat{X}=-X$ and denote
by $\hat{\mathbf{P}}$ for its law starting from $0$.  Note that the
process $(X,\hat{\mathbf{P}})$ is a stable L\'evy process with no
negative jumps of index $\alpha\in (1,2)$ starting from $0$. From
the Lamperti representation, it is clear that $T=I^{*}$ and from the
self-similar property, we have
\[
\begin{split}
\mathbf{P}_1(T>t)&=\mathbf{P}_0(\sigma^{-}_{-1}>t)=\mathbf{P}_0
\left(\inf_{0\leq s\leq t}X_s>-1\right)\\
&=\mathbf{P}_0\left(t^{1/\alpha}\inf_{0\leq s\leq 1}X_s>-1\right)=
\hat{\mathbf{P}}\left(t^{1/\alpha}\sup_{0\leq s\leq 1}X_s<1\right)\\
&=\hat{\mathbf{P}}\left(\left(\frac{1}{\sup_{0\leq s\leq
1}X_s}\right)^{\alpha}> t \right).
\end{split}
\]
Hence the exponential functional $I^{*}$, under $\mathbf{P}_1$, has
the same law as $(\sup_{0\leq s\leq 1}X_{1})^{-\alpha}$, under
 $\hat{\mathbf{P}}$. \\
Recently,   Bernyk, Dalang and Peskir \cite{BDP} computed the
density of the supremum of a stable L\'evy process with no negative
jumps of index $\alpha\in (1,2)$. More precisely with our notation,
the density $f$ of  $\sup_{0\leq s\leq 1}X_{1}$, under
$\hat{\mathbf{P}}$, is described as follows
\[
f(x)=c^{1/\alpha}\sum_{n=1}^{\infty}\frac{\alpha n-1}{\Gamma(\alpha
n)\Gamma(-n+1+1/\alpha)}x^{n\alpha-2},\quad \textrm{for}\quad x>0.
\]
Therefore the density $p^{*}$ of $I^{*}$ is given by
\[
p^{*}(x)=c^{1/\alpha}\sum_{n=1}^{\infty}\frac{\alpha
n-1}{\Gamma(\alpha
n)\Gamma(-n+1+1/\alpha)}x^{\alpha(2-n\alpha)},\quad
\textrm{for}\quad x>0,\] which completes the proof.\QED

\end{document}